\magnification 1000
\baselineskip .85truecm
\hsize 15truecm
\vsize 24truecm
\input amssym.tex


\def \za{\alpha}
\def \zb{\beta}
\def \zg{\gamma}

\def \zl{\lambda}

\def \zp{\pi}

\def \zr{\rho}

\def \zf{\varphi}

\def \zF{\Phi}


\def \zjma{\jmath}
\def \zlma{\ell}


\def \zin{\cap}

\def \zdi{\oplus}

\def \zmm{\pm}

\def \zpor{\times}
\def \zci{\circ}


\def \zmai{\geq}
\def \zco{\subset}

\def \zpe{\in}



\def \zfl{\rightarrow}

\def \zbv{\mid}

\def \z/{\over}

{\bf COMPLEX POLYNOMIAL REPRESENTATION OF $\zp_{n+1}(S^{n})$ AND
$\zp_{n+2}(S^{n})$.}\bigskip

\noindent Francisco-Javier TURIEL

\noindent Geometr'a y Topolog'a, Facultad de Ciencias,
Campus de Teatinos, 29071 M\'alaga, Spain

\noindent e-mail: turiel@agt.cie.uma.es
\bigskip
----------------------------------------------------------------

{\bf Abstract.} The complex affine quadric $Q^{m}=\{ z\zpe {\Bbb
C}^{m+1}\zbv z_{1}^{2}+...+z_{m+1}^{2}=1\}$ deforms by retraction
onto $S^{m}$; this allows us to identify $[Q^{k},Q^{n}]$ and
$[S^{k},S^{n}]=\zp_{k}(S^{n})$. Thus one will say that an element
of $\zp_{k}(S^{n})$ is complex representable if there exists a
complex polynomial map from $Q^{k}$ to $Q^{n}$ corresponding to
this class. In this Note we show that $\zp_{n+1}(S^{n})$ and
$\zp_{n+2}(S^{n})$ are complex representable.

---------------------------------------------------------------
\bigskip

A map from $X\zco {\Bbb K}^{m}$ to  $Y\zco {\Bbb K}^{r}$, ${\Bbb K}={\Bbb R}$ or
${\Bbb C}$, is called {\it polynomial} if
it is the restriction to $X$ of a polynomial (or algebraic) map
$f:{\Bbb K}^{m}\zfl {\Bbb K}^{r}$ such that $f(X)\zco Y$.
Representing elements of each homotopy group $\zp_{k}(S^{n})$, where
$S^{n}=\{ x\zpe {\Bbb R}^{n+1}\zbv x_{1}^{2}+...+x_{n+1}^{2}=1\}$ is the $n$-dimensional sphere,
by real polynomial maps is a classic question now. A first result by P.F. Baum (see[1])
suggested a wide affirmative answer but, later on, R. Wood showed that if $k$ is a power of 2
then all polynomial maps from $S^k$ to $S^{k-1}$ are constant. Consequently the real polynomial representation
is impossible in many cases among them, for example, the non-trivial element of $\zp_{4}(S^{3})$.
Therefore it seems reasonable to consider the complex representation, where this obstacle disappears and
nothing changes up to homotopy.

Let $Q^{m-1}$ be the affine quadric in ${\Bbb C}^{m}$ (often called the complex
sphere) defined by the equation
$z_{1}^{2}+...+z_{m}^{2}=1$. If one set $z_{\zjma}=x_{\zjma}+iy_{\zjma}$, that is
to say, if one identifies ${\Bbb C}^{m}$ to ${\Bbb R}^{2m}={\Bbb R}^{m}\zpor {\Bbb R}^{m}$
endowed with coordinates $(x,y)=(x_{1},...,x_{m},y_{1},...,y_{m})$ then $Q^{m-1}$ is given
by the equations $x_{1}^{2}+...+x_{m}^{2}=y_{1}^{2}+...+y_{m}^{2}+1$, $x_{1}y_{1}+...+x_{m}y_{m}=0$.
The map $g(x,y)=((y_{1}^{2}+...+y_{m}^{2}+1)^{-1/2}x, y)$ gives rise to a
diffeomorphism between $Q^{m-1}$ and the tangent
bundle of the sphere $TS^{m-1}\zco {\Bbb R}^{m}\zpor {\Bbb R}^{m}$. On the other hand
$S^{m-1}= Q^{m-1}\zin ({\Bbb R}^{m}\zpor \{0\})$, and the homotopy
$H(x,y,t)=(((1-t)^{2}y_{1}^{2}+...+(1-t)^{2}y_{m}^{2}+1)^{1/2}(y_{1}^{2}+...+y_{m}^{2}+1)^{-1/2}x,
(1-t)y)$  shows that $Q^{m-1}$ deforms by retraction onto
$S^{m-1}$, which translates the usual retraction of a vector bundle onto the zero section.
Therefore $[Q^{k},Q^{n}]$ is naturally isomorphic to $[S^{k},S^{n}]=\zp_{k}(S^{n})$
under the inclusion of $S^k$ in $Q^k$
and the retraction $H_{1}:Q^{n}\zfl S^{n}$ (see [5] for details). {\it From now on the homotopy
class of a map from $Q^{k}$ to $Q^{n}$ will be regarded as an element of $\zp_{k}(S^{n})$.} We will say
that an element of $\zp_{k}(S^{n})$ is {\it complex representable} if there exists a complex polynomial
map from $Q^k$ to $Q^n$ corresponding to this class; $\zp_{k}(S^{n})$ {\it complex
reprentable} will mean that all its elements are complex representable.
Note that if $f:{\Bbb R}^{m}\zfl {\Bbb R}^{r}$ is a polynomial map and
$f(S^{m-1})\zco S^{r-1}$ then $f(Q^{m-1})\zco Q^{r-1}$ when $f$ is extended to a
polynomial map from ${\Bbb C}^m$ to ${\Bbb C}^r$. In short, real representable implies
complex representable.

Few things are known about the complex representation. In [5] it is proved that $\zp_{n}(S^{n})$
and $\zp_{2k+1}(S^{2k})$ are complex representable. In this Note one shows that $\zp_{n+1}(S^{n})$ and
$\zp_{n+2}(S^{n})$ are complex representable for any dimension $n$.

On each ${\Bbb C}^m$ we consider the symmetric bilinear form $b$ defined
by $b(z,u)=z_{1}u_{1}+...+z_{m}u_{m}$ and its associated quadratic form
$q(z)=b(z,z)$; then $Q^{m-1}=q^{-1}(1)$. A map $f:{\Bbb C}^{m}\zfl {\Bbb C}^{r}$
will called {\it pseudo-homogeneous
of order a natural number $k$} when $q(f(z))=q(z)^{k}$ whatever $z\zpe{\Bbb C}^m$;
in this case $f(Q^{m-1})\zco Q^{r-1}$. Note that order and degree
may be quite different. A polynomial
map from ${\Bbb C}^m$ to ${\Bbb C}^r$, homogeneous of degree $k$ and sending $f(Q^{m-1})$ onto
$Q^{r-1}$, is pseudo-homogeneous of order $k$ as well, but the converse does not hold
unless the polynomial map is real.

{\bf Lemma 1.} {\it For every natural number $k\zmai 1$ there exist three polynomials in two
variables $\zr$, $\zb_1$ and $\zb_2$, the first and the second ones with real coefficients and
the third one with imaginary coefficients,
such that $(t-s)\zr^{2}(s,t)+s^{2k-1}=t^{k}(\zb_{1}^{2}(s,t)+\zb_{2}^{2}(s,t))$ for any
$(s,t)\zpe {\Bbb C}^{2}$, and $\zr(1,t)>0$ when $t\zpe
{\Bbb R}$ et $t\zmai 0$.}

{\bf Proof.} by Lemma 2 of [3], for each natural number $\zlma\zmai 0$, there exist two real
polynomials in one variable $\zf_{\zlma}$
and $\zl_{\zlma}$ such that
$(t-1)\zf_{\zlma}^{2}(t)+1=t^{\zlma+1}\zl_{\zlma}(t)$, the degree of $\zl_{\zlma}$ equals $\zlma$
and $\zf_{\zlma}(t)>0$ whatsoever $t\zmai 0$. So the degree
of $\zf_{\zlma}$ equals $\zlma$ too. Now set $k=\zlma +1$, $\zr(s,t)=s^{\zlma}\zf_{\zlma}(t/s)$,
$\zb(s,t)=s^{\zlma}\zl_{\zlma}(t/s)$, $\zb_{1}=\zb+(1/4)$ and
$\zb_{2}=i\zb-(i/4)$. $\square$

Let $f,g:{\Bbb C}^{m}\zfl {\Bbb C}^{r}$ be two polynomial maps pseudo-homogeneous of order $k\zmai 1$
and $b$-orthogonal (that is to say $b(f(z),g(z))=0$ for each $z\zpe{\Bbb C}^m$).
Then $f,g:Q^{m-1}\zfl Q^{r-1}$ are homotopic by means of the homotopy
${\tilde H}(z,t)=cos(t\zp/2)f(z)+ sin(t\zp/2)g(z)$ among others. Moreover, for any natural number
$\zlma$, we may define a polynomial map
$F:{\Bbb C}^{m+\zlma}={\Bbb C}^{m}\zpor{\Bbb C}^{\zlma}\zfl
{\Bbb C}^{r+\zlma}={\Bbb C}^{r}\zpor{\Bbb C}^{\zlma}$ by setting:

$F(z,u)=(\zb_{1}(z_{1}^{2}+...+z_{m}^{2}+u_{1}^{2}+...+u_{\zlma}^{2}, z_{1}^{2}+...+z_{m}^{2})f(z)
+\zb_{2}(z_{1}^{2}+...+z_{m}^{2}+u_{1}^{2}+...+u_{\zlma}^{2}, z_{1}^{2}+...+z_{m}^{2})g(z),
\zr(z_{1}^{2}+...+z_{m}^{2}+u_{1}^{2}+...+u_{\zlma}^{2}, z_{1}^{2}+...+z_{m}^{2})u)$

\noindent where $\zr$, $\zb_1$ and $\zb_2$ are
like in  Lemma 1 and $(z,u)=(z_{1},...,z_{m},u_{1},...,u_{\zlma})$.

A straightforward calculation shows that $F$ is pseudo-homogeneous of order $2k-1$.

{\bf Theorem 1.} {\it The element of $\zp_{m+\zlma-1}(S^{r+\zlma-1})$ represented by
$F:Q^{m+\zlma-1}\zfl Q^{r+\zlma-1}$ is the (iterated) suspension of the element
of $\zp_{m-1}(S^{r-1})$ represented by $f:Q^{m-1}\zfl Q^{r-1}$.}

{\bf Proof.} Recall that the retraction of $Q^{n-1}$ onto
$S^{n-1}$ was given by
$H_{1}(x,y)=((y_{1}^{2}+...+y_{n}^{2}+1)^{-1/2}x,0)$. To calculate
the class of $F:Q^{m+\zlma-1}\zfl Q^{r+\zlma-1}$ it is enough to
consider $H_{1}\zci F:S^{m+\zlma -1}\zfl S^{r+\zlma -1}$ when
$n=m+\zlma$. As usual ${\Bbb K}^{n_{1}}$ is identified to the
vector subspace ${\Bbb K}^{n_{1}}\zpor \{0\}$ of ${\Bbb
K}^{n_{1}+n_{2}}$, therefore $S^{n_{1}-1}$ and $Q^{n_{1}-1}$ are
subspaces of $S^{n_{1}+n_{2}-1}$ and $Q^{n_{1}+n_{2}-1}$
respectively.

The last part of $H_{1}\zci F_{\zbv S^{m+\zlma-1}}$ equals
$\zr(1, z_{1}^{2}+...+z_{m}^{2})u$ multiplied by a real positive
function, where $(z,u)$ is real as well. But $\zr(1,
z_{1}^{2}+...+z_{m}^{2})>0$ since $z_{1}^{2}+...+z_{m}^{2}$ is
real non-negative, therefore $H_{1}\zci F:S^{m+\zlma-1}\zfl
S^{r+\zlma-1}$ sends equator onto equator, north hemisphere onto
north hemisphere and south one onto south one. This implies that
its homotopy class is the suspension of the class of $H_{1}\zci
F:S^{m+\zlma-2}\zfl S^{r+\zlma-2}$. But, if $\zlma\zmai 2$, this
last map sends equator to equator, north to north and south to
south and we can start the process again. In short the class of
$H_{1}\zci F:S^{m+\zlma-1}\zfl S^{r+\zlma-1}$ is the suspension of
the one of $H_{1}\zci \zf:S^{m-1}\zfl S^{r-1}$, where
$\zf:Q^{m-1}\zfl Q^{r-1}$ is defined by $\zf(z) =
\zb_{1}(1,1)f(z)+\zb_{2}(1,1)g(z)$; note that
$\zb_{1}^{2}(1,1)+\zb_{2}^{2}(1,1)=1$ as it follows from Lemma 1.
For finishing it is enough to show that $f,\zf:Q^{m-1}\zfl
Q^{r-1}$ are homotopic, which is obvious because
$(\zb_{1}(1,1),\zb_{2}(1,1))$ and $(1,0)$ belong to $Q^1$ and this
last  space is path-connected. $\square$

Let us show that $\zp_{n+1}(S^{n})$ and $\zp_{n+2}(S^{n})$ are
complex representable. As the unit element of a homotopy group is
representable by a constant map and
$\zp_{2}(S^{1})=\zp_{3}(S^{1})=0$ we will assume $n\zmai 2$. The
Hopf-Whitehead method, applied to complex multiplication, gives us
the polynomial map $f(x)=(x_{1}^{2}+x_{2}^{2}-x_{3}^{2}-x_{4}^{2},
2x_{1}x_{3}- 2x_{2}x_{4},2x_{1}x_{4}+2x_{2}x_{3})$ from ${\Bbb
R}^4$ to ${\Bbb R}^3$, sending $S^3$ to $S^2$ and whose homotopy
class spans $\zp_{3}(S^{2})={\Bbb Z}$ (Hopf map). Composing on the
right with the real polynomial maps from $S^3$ to itself,
constructed by R. Wood for any topological degree (see theorem 1
of [4]), yields a real representation of $\zp_{3}(S^{2})$; so this
group is complex representable too.

On the other hand the extension of $f$ to ${\Bbb C}^4$, which is
written  $f(z)=(z_{1}^{2}+z_{2}^{2}-z_{3}^{2}-z_{4}^{2},
2z_{1}z_{3}- 2z_{2}z_{4},2z_{1}z_{4}+2z_{2}z_{3})$, sends $Q^3$
onto $Q^2$ and is quadratic homogeneous, therefore it is
pseudo-homogeneous of order 2. But
$g(z)=(2z_{1}z_{4}-2z_{2}z_{3},2z_{1}z_{2}+2z_{3}z_{4},
z_{2}^{2}+z_{4}^{2}-z_{1}^{2}-z_{3}^{2})$ is pseudo-homogeneous of
order 2 as well and $b$-orthogonal to $f$. Now from Theorem 1
follows that $\zp_{n+1}(S^{n})={\Bbb Z}_{2}$, $n\zmai 3$, is
complex representable since the suspension of the class of
$f:S^{3}\zfl S^{2}$ is the non-trivial element of this group.

One has just shown that the non-trivial element of
$\zp_{4}(S^{3})$ may be represented by a complex polynomial map
$\zf:{\Bbb C}^{5}\zfl{\Bbb C}^{4}$ pseudo-homogeneous of order 3.
So $f_{1}=f\zci\zf:{\Bbb C}^{5}\zfl{\Bbb C}^{3}$ and
$g_{1}=g\zci\zf:{\Bbb C}^{5}\zfl{\Bbb C}^{3}$ are
pseudo-homogeneous of order 6; moreover they are $b$-orthogonal
since $f$ and $g$ were $b$-orthogonal. As
$f_{*}:\zp_{4}(S^{3})\zfl\zp_{4}(S^{2})$ is an isomorphism,
$f_{1}$ represents the non-trivial element of $\zp_{4}(S^{2})$.
But the suspension of this element is the non-trivial element of
$\zp_{n+2}(S^{n})$, $n\zmai 3$ (recall that the image of
$S:\zp_{4}(S^{2})\zfl\zp_{5}(S^{3})$ is the kernel of the Hopf
morphism $H:\zp_{5}(S^{3})={\Bbb Z}_{2}\zfl\zp_{3}(S^{3})={\Bbb
Z}$, so $S$ is an isomorphism). Therefore from Theorem 1 follows
that each $\zp_{n+2}(S^{n})$ is complex representable. In short:

{\bf Theorem 2.}{\it The homotopy groups $\zp_{n+1}(S^{n})$ and
$\zp_{n+2}(S^{n})$ are complex representable.}

{\bf Remark.} a) The non-trivial element of $\zp_{5}(S^{3})$ can
be represented by a complex polynomial map $\zF:{\Bbb
C}^{6}\zfl{\Bbb C}^{4}$ pseudo-homogeneous of order some $k$. So
$f_{2}=f\zci\zF:{\Bbb C}^{6}\zfl{\Bbb C}^{3}$ is
pseudo-homogeneous of order $2k$ and represents the non-trivial
element of $\zp_{5}(S^{2})={\Bbb Z}_{2}$ because
$f_{*}:\zp_{5}(S^{3})\zfl\zp_{5}(S^{2})$ is an isomorphism. On the
other hand $g_{2}=g\zci\zF$ is also pseudo-homogeneous of order
$2k$ and $b$-orthogonal to $f_2$ so, by Theorem 1, the suspension
of this element is always complex representable.

A routine calculation shows that this suspension never vanishes.
In other words the element $\za\zpe\zp_{n+3}(S^{n})-\{0\}$,
$n\zmai 2$, such that $2\za =0$ is complex representable because
$\zp_{6}(S^{3})={\Bbb Z}_{12}$, $\zp_{7}(S^{4})={\Bbb Z}\zdi {\Bbb
Z}_{12}$ and $\zp_{n+3}(S^{n})={\Bbb Z}_{24}$, $n\zmai 5$.

b) Since $\zp_{1}(S^{1})$ may be represented by real homogeneous
polynomial maps and each of them has an orthogonal map of the same
kind (take $g=(-f_{2},f_{1})$ when  $f=(f_{1},f_{2})$), by
complexifying we can apply Theorem 1; therefore $\zp_{n}(S^{n})$
is complex representable by means of pseudo-homogeneous maps of
order odd. More exactly, if the absolute value of the topological
degree of an element of $\zp_{n}(S^{n})-\{0\}$ equals $k$, then
one can take $2k-1$ as order (when $n=2$ any topological degree
$\zmm k$ can be represented by
a complex polynomial map of degree $2k-1$, which is not
pseudo-homogeneous; see[2]).

Note that the odd order is essential for even dimension. Indeed,
if $\zf:{\Bbb C}^{n+1}\zfl{\Bbb C}^{n+1}$ is pseudo-homogeneous of
order $2r$ then $H(z,t)=\zg(t)^{-2r}\zf(\zg(t)z)$, where $\zg(t)$
is a path from $1$ to $-1$ in ${\Bbb C}-\{0\}$, is a homotopy
between $\zf:Q^{n}\zfl Q^{n}$ and $\zf\zci (-Id):Q^{n}\zfl Q^{n}$;
consequently the class of $\zf$ vanishes for $n$ even.

c) The Hopf maps from  $S^7$ and $S^{15}$ to $S^4$ and $S^8$
respectively have no $b$-orthogonal map of the same kind. Indeed,
if not the element $\za\zpe\zp_{8}(S^{5})$ corresponding to the
suspension of the Hopf map from $S^7$ to $S^4$ (the other case is
analogous) can be represented by a pseudo-homogeneous map
$\zf:Q^{8}\zfl Q^{5}$ of order 3, which is homotopic to $-\zf\zci
(-Id) :Q^{8}\zfl Q^{5}$ by considering
$H(z,t)=\zg(t)^{-3}\zf(\zg(t)z)$ and the path $\zg(t)$ as
previously. Therefore $2\za=0$, {\it contradiction} because $\za$
spans $\zp_{8}(S^{5})$.

{\bf References}
\bigskip

\noindent 1.  P. F. Baum, {\it Quadratic maps and stable homotopy groups of spheres}, Illinois J.Math.
{\bf 11} (1967), 586--595.

\noindent 2.  G.M. Golasi{\'n}ski and F. G\'omez Ruiz, {\it Polynomial and regular maps into
Grassmannians}, K-Theory {\bf 26} (2002), 51-58.

\noindent 3.  F.J. Turiel, {\it Polynomial maps and even dimensional spheres}, Proc.
Amer. Math.Soc. (2007) to appear.

\noindent 4.  R. Wood, {\it Polynomial maps from spheres to spheres},
Invent. Math. {\bf 5} (1968), 163--168.

\noindent 5.  R. Wood, {\it Polynomial maps of affine quadrics},
Bull. London Math. Soc. {\bf 25} (1993), 491--497.



\end